\documentclass[twoside,11pt]{article}

\RequirePackage[colorlinks=true]{hyperref}
\usepackage{amsmath,amsfonts,bm,amsthm}
\usepackage[capitalize,noabbrev]{cleveref}

\usepackage[abbrvbib,nohyperref]{jmlr2e}

\usepackage[dvipsnames]{xcolor}
\usepackage{inputenc}
\usepackage{mathtools}
\usepackage{bbm}
\usepackage{booktabs}
\usepackage{enumerate}
\usepackage{float}
\usepackage{algorithm}
\usepackage[noEnd=true]{algpseudocodex}
\usepackage{microtype}
\usepackage{enumitem}
\usepackage{physics}
\algrenewcommand{\algorithmicrequire}{\textbf{Input:}}
\algrenewcommand{\algorithmicensure}{\textbf{Output:}}
\usepackage[most]{tcolorbox}
\tcbset{on line,
  boxsep=2pt, left=0pt, right=0pt, top=1pt, bottom=1pt,
  colframe=black, colback=BurntOrange!80,
  leftrule=1pt, rightrule=1pt, toprule=1pt, bottomrule=1pt,
  highlight math style={enhanced}
}
\usepackage{textcomp}

\definecolor{blu}{RGB}{50,106,170}
\definecolor{gre}{RGB}{63,129,71}
\definecolor{drd}{RGB}{179,45,38}
\definecolor{pur}{RGB}{81,64,130}
\hypersetup{hidelinks}

\newcommand{\diff}{\,\text{d}}
\newcommand{\T}{\mathsf{T}}
\newcommand{\lsqrt}[1]{\sqrt{#1}}

\newcommand{\tria}{\operatorname{tria}}

\newcommand{\rone}[1]{{#1}}
\newcommand{\rtwo}[1]{{#1}}

\usepackage{pifont}

\usepackage{lastpage}
\jmlrheading{25}{2024}{1-\pageref{LastPage}}{10/23; Revised 6/24}{7/24}{23-1261}{Nathanael Bosch, Adrien Corenflos, Fatemeh Yaghoobi, Filip Tronarp, Philipp Hennig and Simo S\"{a}rkk\"{a}}

\ShortHeadings{%
  Parallel-in-Time Probabilistic Numerical ODE Solvers%
}{%
  Bosch, Corenflos, Yaghoobi, Tronarp, Hennig, Särkkä%
}

\title{Parallel-in-Time Probabilistic Numerical ODE Solvers}
\author{%
  \name Nathanael Bosch \email nathanael.bosch@uni-tuebingen.de \\
  \addr Tübingen AI Center, University of Tübingen
  \AND
  \name Adrien Corenflos \email adrien.corenflos.stats@gmail.com \\
  \addr Department of Electrical Engineering and Automation, Aalto University
  \AND
  \name Fatemeh Yaghoobi \email fatemeh.yaghoobi@aalto.fi \\
  \addr Department of Electrical Engineering and Automation, Aalto University
  \AND
  \name Filip Tronarp \email filip.tronarp@matstat.lu.se \\
  \addr Center for Mathematical Sciences, Lund University
  \AND
  \name Philipp Hennig \email philipp.hennig@uni-tuebingen.de \\
  \addr Tübingen AI Center, University of Tübingen
  \AND
  \name Simo Särkkä \email simo.sarkka@aalto.fi \\
  \addr Department of Electrical Engineering and Automation, Aalto University
}
\editor{Ryan Adams}
\firstpageno{1}

\begin{document}
\maketitle

\begin{abstract}%
  Probabilistic numerical solvers for ordinary differential equations (ODEs) treat the numerical simulation of dynamical systems as problems of Bayesian state estimation.
  Aside from producing posterior distributions over ODE solutions and thereby quantifying the numerical approximation error of the method itself, one less-often noted advantage of this formalism is the algorithmic flexibility gained by formulating numerical simulation in the framework of Bayesian filtering and smoothing.
  In this paper, we leverage this flexibility and build on the time-parallel formulation of iterated extended Kalman smoothers to formulate a \emph{parallel-in-time} probabilistic numerical ODE solver.
  Instead of simulating the dynamical system sequentially in time, as done by current probabilistic solvers, the proposed method processes all time steps in parallel and thereby reduces the \rone{computational complexity} from \emph{linear} to \emph{logarithmic} in the number of time steps.
  We demonstrate the effectiveness of our approach on a variety of ODEs and compare it to a range of both classic and probabilistic numerical ODE solvers.
\end{abstract}


\section{Introduction}

Ordinary differential equations (ODEs) are used throughout the sciences to describe the evolution of dynamical systems over time.
In machine learning, ODEs provide a continuous description of certain neural networks \citep{chen2018} and optimization procedures \citep{helmke2012optimization,su2016},
and are used in generative modeling with normalizing flows \citep{papamakarios2021} and diffusion models \citep{song2021scorebased}, among others.
Unfortunately, all but the simplest ODEs are too complex to be solved analytically. Therefore, numerical methods are required to obtain a solution.
While a multitude of numerical solvers has been developed over the last century \citep{hairer2008solving,deuflhardbornemann,butcher2016numerical}, most commonly-used methods do not provide a quantification of their own inevitable numerical approximation error.

Probabilistic numerics provides a framework for treating classic numerical problems as problems of probabilistic inference
\citep{hennig15_probab_numer_uncer_comput,oates2019,hennig_osborne_kersting_2022}.
In the context of ODEs, methods based on Gaussian process regression
\citep{skilling1992,hennig14}
and in particular Gauss--Markov regression
\citep{Schober2019,Kersting2020a,Tronarp2019c}
provide an efficient and flexible approach to compute posterior distributions over the solution of ODEs
\citep{Bosch2021,Kramer2020stable},
and even partial differential equations \citep{kramer2022pde,bosch2023probabilistic}
and differential-algebraic equations \citep{Bosch2022}.
These so-called \emph{ODE filters} typically scale cubically in the ODE dimension (as do most \emph{implicit} ODE solvers) and specific approximations enable linear scaling (shared by most \emph{explicit} solvers) \citep{Kramer2022highdim}.
But to date, their linear scaling with the number of time steps remains.

For very large-scale simulations with very long time horizons, the sequential processing in time of most ODE solvers can become a bottleneck.
This motivates the development of \emph{parallel-in-time} methods:
By leveraging the ever-increasing parallelization capabilities of modern computer hardware, parallel-in-time methods can achieve \emph{sub-linear} scaling in the number of time steps \citep{gander-50y}.
One well-known method of this kind is Parareal \citep{LIONS2001}.
It achieves temporal parallelism by combining an expensive, accurate solver with a cheap, coarse solver, in such a way that the fine solver is only ever applied to individual time slices in a parallel manner, leading to
a square-root scaling (in ideal conditions).
But, due to its sequential coarse-grid solve, Parareal still has only limited concurrency
\citep{Gander2007},
and while it has recently been extended probabilistically by \citet{Pentland2021,Pentland2022} to improve its performance and convergence,
these methods do not provide probabilistic solutions to ODEs per se.

In this paper, we leverage the time-parallel formulation of Gaussian filters and smoothers \citep{Sarkka2021,Yaghoobi2021,Yaghoobi2022} and develop a \emph{parallel-in-time} probabilistic numerical ODE solver.
The paper is structured as follows.
\Cref{sec:background} formulates numerical ODE solutions as Bayesian state estimation problems and presents the established, sequential, filtering-based probabilistic ODE solvers.
\Cref{sec:method} then
presents our proposed parallel-in-time method;
first as exact inference for affine ODEs, then as an iterative, approximate algorithm for general nonlinear ODEs.
\Cref{sec:experiments} then presents experiments on a variety of ODEs and compares the performance of our proposed method to that of existing \rone{(sequential)} probabilistic and non-probabilistic ODE solvers.
Finally, \cref{sec:conclusion} concludes with a discussion of our results and an outlook on future work.

\section{Numerical ODE Solutions as Bayesian State Estimation}
\label{sec:background}
Consider an initial value problem (IVP) \rtwo{given by a first-order ODE} of the form
  \begin{align}
    \label{eq:ivp}
    \dot{y}(t) = f(y(t), t), \quad t \in [0, T], \qquad y(0) = y_0,
  \end{align}
with vector field \(f: \mathbb{R}^d \times \mathbb{R} \to \mathbb{R}^d\) and initial value \(y_0 \in \mathbb{R}^d\).
To capture the numerical error that arises from temporal discretization, the quantity of interest in probabilistic numerics for ODEs is the \emph{probabilistic numerical ODE solution}, defined as
\begin{equation}
  \label{eq:pnsol}
  p \left( y(t) ~\Big|~ y(0) = y_0, \left\{ \dot{y}(t_n) = f(y(t_n), t_n) \right\}_{n=1}^N \right),
\end{equation}
for some prior \(p\left(y(t)\right)\) and with \( \{t_n\}_{n=1}^N \subset [0, T]\) the chosen time-discretization.

In the following, we pose the probabilistic numerical ODE solution as a problem of Bayesian state estimation, and we define the prior, likelihood, data, and approximate inference scheme.
For a more detailed description of the transformation of an IVP into a Gauss--Markov regression problem, refer to \citet{Tronarp2019c}.

\subsection{Gauss--Markov Process Prior}
\label{sec:prior}
We model the solution \(y\) of the IVP with a \(\nu\)-times integrated Wiener process prior (IWP(\(\nu\))).
More precisely, let
\(Y(t) = \left[ Y^{(0)}(t), Y^{(1)}(t), \dots, Y^{(\nu)}(t) \right]\)
be the solution of the following linear, time-invariant stochastic differential equation with Gaussian initial condition
\begin{subequations}
  \begin{align}
    \diff Y^{(i)}(t) &= Y^{(i+1)}(t) \diff t, \qquad i = 0, \dots, \nu-1, \\
    \diff Y^{(\nu)}(t) &= \rtwo{\sigma \diff W(t),} \\
    Y(0) &\sim \mathcal{N} \left( \mu_0, \Sigma_0 \right),
  \end{align}
\end{subequations}
with initial mean and covariance \(\mu_0 \in \mathbb{R}^{d(\nu+1)}\), \(\Sigma_0 \in \mathbb{R}^{d(\nu+1) \times d(\nu+1)}\),
\rtwo{diffusion coefficient \(\sigma \in \mathbb{R}_+\)}, and
\(d\)-dimensional Wiener process
\(W: \mathbb{R} \to \mathbb{R}^d\).
Then, \(Y^{(i)}\) is chosen to model the \(i\)-th derivative of the IVP solution \(y\).
By construction, accessing the \(i\)-th derivative can be done by multiplying the state \(Y\) with a projection matrix \(E_i \coloneqq I_d \otimes e_i\), that is, \(Y^{(i)}(t) = E_i Y(t)\).

This continuous-time prior satisfies discrete transition densities
\citep{sarkka_solin_2019}
\begin{equation}
  Y(t+h) \mid Y(t) \sim \rtwo{\mathcal{N} \left( \Phi(h) Y(t), \sigma^2 Q(h) \right)},
\end{equation}
with transition matrix and process noise covariance
\(\Phi(h), Q(h) \in \mathbb{R}^{d(\nu+1) \times d(\nu+1)}\) and step \(h \in \mathbb{R}_+\).
For the IWP(\(\nu\)) these can be computed in closed form
\citep{Kersting2020a}, as
\begin{subequations}
  \begin{align}
    \Phi(h) &= I_d \otimes \breve{\Phi}(h),
    &
      \left[ \breve{\Phi}(h) \right]_{ij} &= \mathbbm{1}_{i=j} \frac{h^{i-j}}{(j-i)!},
    \\
    Q(h) &= I_d \otimes \breve{Q}(h),
    &
    \left[ \breve{Q}(h) \right]_{ij} &= \rone{\frac{h^{2 \nu + 1 - i - j}}{(2 \nu + 1 -i -j)(\nu-i)!(\nu-j)!}.}
  \end{align}
\end{subequations}

\begin{remark}[Alternative Gauss--Markov priors]
  While \(\nu\)-times integrated Wiener process priors have been the most common choice for filtering-based probabilistic ODE solvers in recent years, the methodology is not limited to this choice.
  Alternatives include the \(\nu\)-times integrated Ornstein--Uhlenbeck process
  and the class of Mat\'ern processes, both of which have a similar continuous-time SDE representation as well as Gaussian transition densities in discrete time.
  Refer to \citet{bosch2023probabilistic,Tronarp2021a,sarkka_solin_2019}.
\end{remark}

The initial distribution \(\mathcal{N}(\mu_0, \Sigma_0)\) is chosen such that it encodes the initial condition \(y(0) = y_0\).
Furthermore, to improve the numerical stability and the quality of the posterior,
we initialize not only on the function value
\(Y^{(0)}(0) = y_0\),
but also the higher order derivatives, that is,
\(Y^{(i)}(0) = \frac{\diff^i y}{\diff t^i}(0)\) for all \(i \leq \nu\)
\citep{Kramer2020stable}.
These terms can be efficiently computed via Taylor-mode automatic differentiation
\citep{griewank2000evaluating,bettencourt2019taylormode}.
As a result, we obtain an initial distribution with mean
\begin{equation}
  \label{eq:mu0}
  \mu_0 = \left[ y_0, \frac{\diff y}{\diff t}(0), \dots, \frac{\diff^{\nu} y}{\diff t^{\nu}}(0) \right]^T,
\end{equation}
and zero covariance \(\Sigma_0 = 0\), since the initial condition has to hold exactly.

\subsection{Observation Model and Data}
To relate the introduced Gauss--Markov prior to the IVP problem from
\cref{eq:ivp},
we define an observation model in terms of the information operator
\begin{equation}
  \mathcal{Z}[y](t) \coloneqq \dot{y}(t) - f\left( y(t), t \right).
\end{equation}
By construction, \(\mathcal{Z}\) maps the true IVP solution \(y\) \emph{exactly} to the zero function,
that is, \(\mathcal{Z}[y] \equiv 0\).
In terms of the continuous process \(Y\), the information operator can be expressed as
\begin{equation}
  \mathcal{Z}[Y](t) = E_1 Y(t) - f\left( E_0 Y(t), t \right),
\end{equation}
where \(E_0\) and \(E_1\) are the projection matrices introduced in \cref{sec:prior} which select the zeroth and first derivative from the process \(Y\), respectively.
There again, if \(Y\) corresponds to the true IVP solution (and its true derivatives), then \(\mathcal{Z}[Y] \equiv 0\).

Conversely, inferring the true IVP solution requires conditioning the process \(Y(t)\) on \(Z(t)=0\) over the whole continuous interval \(t \in [0, T]\).
Since this is in general intractable, we instead condition \(Y(t)\) only on discrete observations \(Z(t_n) = 0\) on a grid \( \mathbb{T} = \{t_n\}_{n=1}^N \).
This leads to the Dirac likelihood model commonly used in ODE filtering
\citep{Tronarp2019c}:%
\begin{equation}
  \label{eq:observation-model}
  Z(t_n) \mid Y(t_n) \sim \delta \left( Y^{(1)}(t_n) - f\left( Y^{(0)}(t_n), t_n \right) \right),
\end{equation}
with zero-valued data \(Z(t_n) = 0\) for all \(t_n \in \mathbb{T}\).

\begin{remark}[Information operators for other differential equation problems]
  Similar information operators can be defined for other types of differential equations that are not exactly of the first-order form as given in \cref{eq:ivp}, such as higher-order differential equations, Hamiltonian dynamics, or differential-algebraic equations
  \citep{Bosch2022}.
\end{remark}

\subsection{Discrete-Time Inference Problem}
The combination of prior, likelihood, and data results in a Bayesian state estimation problem%
\begin{subequations}
  \label{eq:inference-problem}
  \begin{align}
    Y(0) &\sim \mathcal{N} \left( \mu_0, \Sigma_0 \right), \\
    Y(t_{n+1}) \mid Y(t_n) &\sim \mathcal{N} \left( \Phi(t_{n+1} - t_n) Y(t_n), Q(t_{n+1} - t_n) \right), \\
    \label{eq:inference-problem-obs}
    Z(t_n) \mid Y(t_n) &\sim \delta \left( Y^{(1)}(t_n) - f\left( Y^{(0)}(t_n), t_n \right) \right),
  \end{align}
\end{subequations}
with zero data
\(Z(t_n) = 0\) for all
\(t_n \in \mathbb{T}\).
The posterior distribution over \(Y^{(0)}(t)\) then provides a probabilistic numerical ODE solution to the given IVP, as formulated in \cref{eq:pnsol}.

This
is a standard nonlinear Gauss--Markov regression problem, for which many approximate inference algorithms have previously been studied
\citep{sarkka2023bfs2}.
In the context of probabilistic ODE solvers,
a popular approach for efficient approximate inference is Gaussian filtering and smoothing, where the solution is approximated with Gaussian distributions
\begin{equation}
  p \left( Y(t) \mid \{Z(t_n) = 0\}_{n=1}^N \right) \approx \mathcal{N} \left( \mu(t), \Sigma(t) \right).
\end{equation}
This is most commonly performed with extended Kalman filtering (EKF) and smoothing
(EKS) \citep{Schober2019, Tronarp2019c, Kersting2020a};
though other methods have been proposed, for example based on numerical quadrature
\citep{kersting16} or particle filtering \citep{Tronarp2019c}.
\emph{Iterated} extended Kalman smoothing
\citep[e.g.][]{Bell1994,sarkka2023bfs2}
computes the ``maximum a posteriori'' estimate of the probabilistic numerical ODE solution
\citep{Tronarp2021a}.
This will be the basis for the parallel-in-time ODE filter proposed in this work, explained in detail in \cref{sec:method}.

\subsection{Practical Considerations for Probabilistic Numerical ODE Solvers}
\label{sec:practical-considerations}
While Bayesian state estimation methods such as the extended Kalman filter and smoother can, in principle, be directly applied to the formulated state estimation problem, there are a number of modifications and practical considerations that should be taken into account:

\begin{itemize}
\item \emph{Square-root formulation:}
  \label{par:sqrt-comment}
  Gaussian filters often suffer from numerical stability issues when applied to the ODE inference problem defined in \cref{eq:inference-problem}, in particular when using high orders and small steps,
  \rone{
  due to the ill-conditioning of the state transition covariance \(Q\)
  and numerical round-off error from finite precision arithmetic
  \citep{Kramer2020stable}.
  }%
  To alleviate these issues, probabilistic numerical ODE solvers are typically formulated in square-root form \citep{Kramer2020stable}; this is also the case for the proposed parallel-in-time method.
\item \emph{Preconditioned state transitions:}
  \citet{Kramer2020stable} suggest a coordinate change preconditioner to make the state transition matrices step-size independent and thereby improve the numerical stability of EKF-based probabilistic ODE solvers.
  This preconditioner is also used in this work.
\item \emph{Uncertainty calibration:}
  The Gauss--Markov prior as introduced in \cref{sec:prior} has a free parameter, the \rtwo{diffusion \(\sigma\)}, which directly influences the uncertainty estimates returned by the ODE filter, but not its mean estimates.
  \rtwo{In this paper,
  we compute a quasi-maximum likelihood estimate for the parameter \(\sigma\) post-hoc, as suggested by \citet{Tronarp2019c}.}
\item \emph{Approximate linearization:}
  Variants of the standard EKF/EKS-based inference have been proposed in which the linearization of the vector-field is done only approximately.
  Approximating the Jacobian of the ODE vector field with \rone{a zero matrix} enables inference with a complexity which scales only linearly with the ODE dimension \citep{Kramer2022highdim}, while still providing polynomial convergence rates
  \citep{Kersting2020a}.
  A diagonal approximation of the Jacobian preserves the linear complexity, but improves the stability properties of the solver
  \citep{Kramer2022highdim}.
  In this work, we only consider the exact first-order Taylor linearization.
\item \emph{Local error estimation and step-size adaptation:}
  Rather than predefining the time discretization grid, certain solvers employ an adaptive approach where the solver dynamically constructs the grid while controlling an internal estimate of the numerical error.
  Step-size adaptation based on \emph{local} error estimates have been proposed for both classic
  \citep[Chapter II.4]{hairer2008solving}
  and probabilistic ODE solvers
  \citep{Schober2019,Bosch2021}.
  On the other hand, \emph{global} step-size selection is often employed in numerical boundary value problem (BVP) solvers \citep[Chapter 9]{ascher1995},
  and has been extended to filtering-based probabilistic BVP solvers \citep{kraemer2021bvp}.
  For our purposes, we will focus on fixed grids.
\end{itemize}

\section{Parallel-in-Time Probabilistic Numerical ODE Solvers}
\label{sec:method}
This section develops the main method proposed in this paper:
a parallel-in-time probabilistic numerical ODE solver.

\subsection{\rone{Time-Parallel} Exact Inference in Affine Vector Fields}
\label{sec:affine}

Let us first consider the simple case: An initial value problem with affine vector field
\begin{align}
  \label{eq:affine-ivp}
  \dot{y}(t) = L(t) y(t) + d(t), \quad t \in [0, T], \qquad
  y(0) = y_0.
\end{align}
The corresponding information model of the probabilistic solver is then also affine, with
\begin{subequations}
\begin{align}
  Z(t) \mid Y(t) &\sim \delta \left( H(t) Y(t) - d(t) \right), \\
  H(t) &\coloneqq E_1 - L(t) E_0.
\end{align}
\end{subequations}
Let \(\mathbb{T} = \{t_n\}_{n=1}^N \subset [0, T]\) be a discrete time grid.
To simplify the notation in the following, we will denote a function evaluated at time \(t_n\) by a subscript \(n\), that is \(Y(t_n) =: Y_n\),
except for the transition matrices where we will use
\(\Phi_n \coloneqq \Phi(t_{n+1} - t_n)\) and \(Q_n \coloneqq Q(t_{n+1} - t_n)\).
Then, the Bayesian state estimation problem from \cref{eq:inference-problem} reduces to inference of \(Y(t)\) in the model
\begin{subequations}
  \label{eq:affine-estimation}
  \begin{align}
    Y_0 &\sim \mathcal{N} \left( \mu_0, \Sigma_0 \right), \\
    Y_{n+1} \mid Y_n &\sim \mathcal{N} \left( \Phi_n Y_n, Q_n \right), \\
    Z_n \mid Y_n &\sim \delta \left( H_n Y_n - d_n \right),
  \end{align}
\end{subequations}
with zero data \(Z_n = 0\) for all \(n = 1, \dots, N\).
Since this is an affine Gaussian state estimation problem, it can be solved exactly with Gaussian filtering and smoothing
\citep{kalman1960,rauchtungstriebel1965,sarkka2023bfs2};
see also \citep{Tronarp2019c,Tronarp2021a} for explicit discussions of probabilistic numerical solvers for affine ODEs.

Recently, \citet{Sarkka2021} presented a \rone{time-parallel} formulation of Bayesian filtering and smoothing,
as well as a concrete algorithm for exact linear Gaussian filtering and smoothing---which could be directly applied to the problem formulation in \cref{eq:affine-estimation}.
But as mentioned in \cref{par:sqrt-comment}, the resulting ODE solver might suffer from numerical instabilities.
Therefore, we use the square-root formulation of the \rone{time-parallel} linear Gaussian filter and smoother by
\citet{Yaghoobi2022}.
In the following, we review the details of the algorithm.

\subsubsection{\rone{Time-Parallel} General Bayesian Filtering and Smoothing}

First, we follow the presentation of \citet{Sarkka2021} and formulate Bayesian filtering and smoothing as prefix sums.
We define elements
\(a_n = (f_n, g_n)\) with
\begin{subequations}
\begin{align}
  f_n(Y_n \mid Y_{n-1}) &= p(Y_n \mid Z_n, Y_{n-1}), \\
  g_n(Y_{n-1}) &= p(Z_n \mid Y_{n-1}),
\end{align}
\end{subequations}
where for \(n=1\) we have \(p(Y_1 \mid Z_1, Y_0) = p(Y_1 \mid Z_1)\) and \(p(Z_1 \mid Y_0) = p(Z_1)\),
together with a binary \rone{filtering} operator
\(\otimes_f : (f_i, g_i) \otimes_f (f_j , g_j ) \mapsto (f_{ij} , g_{ij})\)
defined by
\begin{subequations}
\begin{align}
  f_{ij}(x \mid z) &\coloneqq \frac{\int g_j(y) f_j(x \mid y) f_i(y \mid z) \diff y}{\int g_j(y) f_i(y \mid z) \diff y}, \\
  g_{ij}(z) &\coloneqq g_i(z) \int g_j(y) f_i(y \mid z) \diff y.
\end{align}
\end{subequations}
Then,
\citet[Theorem 3]{Sarkka2021}
show that \(\otimes_f\) is associative and that
\begin{align}
  a_1 \otimes_f \dots \otimes_f a_n = \begin{bmatrix} p(Y_n \mid Z_{1:n}) \\ p(Z_{1:n}) \end{bmatrix},
\end{align}
that is, the filtering marginals and the marginal likelihood of the observations at step \(n\) are the results of a cumulative sum of the elements \(a_{1:n}\) under \(\otimes_f\).
Since the \rone{filtering} operator \(\otimes_f\) is associative, this quantity can be computed in parallel with prefix-sum algorithms, such as the parallel scan algorithm by \citet{Blelloch1989}.

\begin{remark}[On Prefix-Sums]
  Prefix sums, also known as cumulative sums or inclusive scans, play an important role in parallel computing.
  Their computation can be efficiently parallelized and, if enough parallel resources are available, their (span) computational cost can be reduced from \emph{linear} to \emph{logarithmic} in the number of elements.
  One such algorithm is the well-known parallel scan algorithm by \citet{Blelloch1989} which, given \(N\) elements and \(N/2\) processors, computes the prefix sum in \(2 \lceil \log_2 N \rceil\) sequential steps with \(2N-2\) invocations of the binary operation.
  This algorithm is implemented in both tensorflow \citep{tensorflow2015-whitepaper} and JAX \citep{jax2018github};
  the latter is used in this work.
\end{remark}

The time-parallel smoothing step can be constructed similarly:
We define elements
\(b_n = p(Y_n \mid Z_{1:n}, Y_{n+1})\), with \(b_N = p(Y_N \mid Z_{1:N})\),
and a binary \rone{smoothing} operator
\(b_i \otimes_s b_j = b_{ij}\), with
\begin{align}
  b_{ij}(x \mid z) &= \int b_i(x \mid y) b_j(y \mid z) \diff y.
\end{align}
Then, \(\otimes_s\) is associative and the smoothing marginal at time step \(n\) is the result of a reverse cumulative sum of the elements \(b_{n:N}\) under \(\otimes_s\)
\citep{Sarkka2021}:%
\begin{align}
  b_n \otimes_s \dots \otimes_s b_N = p(Y_n \mid Z_{1:N}).
\end{align}
Again, since the smoothing operator \(\otimes_s\) is associative, this cumulative sum can be computed in parallel with a prefix-sum algorithm
\citep{Blelloch1989}.

\subsubsection{\rone{Time-Parallel} Linear Gaussian Filtering in Square-Root Form}

In the linear Gaussian case, the filtering elements
\(a_n = (f_n, g_n)\) can be parameterized by a set of parameters
\(\{A_n, b_n, C_n, \eta_n, J_n\}\) as follows:
\begin{subequations}
\begin{align}
  f_n(Y_n \mid Y_{n-1})
  &= p(Y_n \mid Z_n, Y_{n-1})
  = \mathcal{N} \left( Y_n; A_n Y_{n-1} + b_n, C_n \right), \\
  g_n(Y_{n-1})
  &= p(Z_n \mid Y_{n-1})
  \propto \mathcal{N}_I \left( Y_{n-1};  \eta_n, J_n \right),
\end{align}
\end{subequations}
where \(\mathcal{N}_I\) denotes a Gaussian density parameterized in information form,
that is, \(\mathcal{N}_I(x; \eta, J) = \mathcal{N}(x; J^{-1} \eta, J^{-1})\).
The parameters
\(\{A_n, b_n, C_n, \eta_n, J_n\}\) can be computed explicitly from the given state-space model \citep[Lemma 7]{Sarkka2021}.
But since probabilistic numerical ODE solvers require a numerically stable implementation of the underlying filtering and smoothing algorithm \citep{Kramer2020stable}, we formulate the \rone{time-parallel} linear Gaussian filtering algorithm in square-root form, following \citet{Yaghoobi2022}.

To this end,
let \(\lsqrt{M}\) denote a left square-root of a positive semi-definite matrix \(M\), that is, \(\lsqrt{M} \lsqrt{M}^T = M\);
the matrix \(\lsqrt{M}\) is sometimes also called a ``generalized Cholesky factor'' of \(M\) \citep{grewal2014}.
To operate on square-root matrices, we also define the \emph{triangularization} operator:
Given a wide matrix \(M \in \mathbb{R}^{n \times m}\), \(m \geq n\), the triangularization operator \(\tria(M)\) first computes the QR decomposition of
\(M^\T\), that is, \(M^\T = QR\), with wide orthonormal \(Q \in \mathbb{R}^{m \times n}\) and square upper-triangular \(R \in \mathbb{R}^{n \times n}\), and then returns \(R^\T\).
This operator plays a central role in square-root filtering algorithms as it enables the numerically stable addition of covariance matrices, provided square-roots are available.
Given two positive semi-definite matrices \(A, B \in \mathbb{R}^{n \times n}\) with square-roots \(\lsqrt{A}, \lsqrt{B}\), a square-root of the sum \(A + B\) can be computed as
\begin{align}
  \lsqrt{A + B} = \tria \left( \begin{bmatrix} \lsqrt{A} & \lsqrt{B} \end{bmatrix} \right).
\end{align}
With these definitions in place, we briefly review the \rone{time-parallel} linear Gaussian filtering algorithm in square-root form as provided by \citet{Yaghoobi2022} in the following.

\paragraph{Parameterization of the filtering elements.}
Let \(m_0 = \mu_0\), \(P_0 = \Sigma_0\), and \(m_n = 0\), \(P_n = 0\) for all \(n \geq 1\), and define
\begin{subequations}
\begin{align}
  m_n^- &= \Phi_{n-1} m_{n-1}, \\
  \lsqrt{P_n^-} &= \tria \left( \begin{bmatrix} \Phi_{n-1} \lsqrt{P_{n-1}} & \lsqrt{Q_{n-1}} \end{bmatrix} \right).
\end{align}
\end{subequations}
Then, the square-root parameterization of the filtering elements \(a_n\) is given by
\begin{subequations}
\label{eq:sqrt-filtering-elements}
\begin{align}
  A_n &= (I - K_n H_n) \Phi_{n-1}, \\
  b_n &= m_n^- - K_n \left( H_n m_n^- - d_n \right), \\
  \lsqrt{C_n} &= \Psi_{22}, \\
  \eta_n &= \lsqrt{J_n} \lsqrt{S_n}^{-1} d_n, \\
  \lsqrt{J_n} &= \Phi_{n-1}^\T H_n^\T \sqrt{S_n}^{-\T},
\end{align}
\end{subequations}
where \(I\) is the identity matrix and \(\Psi_{22}\), \(\lsqrt{S_n}\) and \(K_n\) are defined via
\begin{subequations}
\begin{align}
  \begin{bmatrix}
    \Psi_{11} & 0 \\
    \Psi_{21} & \Psi_{22}
  \end{bmatrix}
  &=
  \operatorname{tria} \left(
  \begin{bmatrix}
    H_n \lsqrt{P_n^-} & \lsqrt{R_n} \\
    \lsqrt{P_n^-} & 0
  \end{bmatrix}
  \right), \\
  \lsqrt{S_n} &= \Psi_{11}, \\
  K_n &= \Psi_{21} \Psi_{11}^{-1}.
\end{align}
\end{subequations}
For generality the formula includes an observation noise covariance \(R_n\),
but note that in the context of probabilistic ODE solvers we have a noiseless measurement model with \(\lsqrt{R_n} = 0\).

\paragraph{Associative filtering operator.}
Let \(a_i, a_j\) be two filtering elements, parameterized in square-root form by
\(a_i = \{A_i, b_i, \lsqrt{C_i}, \eta_i, \lsqrt{J_i}\}\) and \(a_j = \{A_j, b_j, \lsqrt{C_j}, \eta_j, \lsqrt{J_j}\}\).
Then, the associative filtering operator \(\otimes_f\)
computes the filtering element \(a_{ij} = a_i \otimes_f a_j\)
as
\begin{subequations}
  \label{eq:sqrt-associative-filtering}
  \begin{align}
  A_{ij} &= A_j A_i - A_j \lsqrt{C_i} \Xi_{11}^{-\T} \Xi_{21}^\T A_i, \\
  b_{ij} &= A_j \left( I - \lsqrt{C_i} \Xi_{11}^{-\T} \Xi_{21}^\T \right) (b_i + \lsqrt{C_i} \lsqrt{C_i}^\T \eta_j) + b_j, \\
  \lsqrt{C_{ij}} &= \tria \left( \begin{bmatrix} A_j \lsqrt{C_i} \Xi_{11}^{-\T} & \lsqrt{C_j} \end{bmatrix} \right), \\
  \eta_{ij} &= A_i^\T \left( I - \Xi_{21} \Xi_{11}^{-1} \lsqrt{C_i}^\T \right) \left( \eta_j - \lsqrt{J_j} \lsqrt{J_j}^\T b_i \right) + \eta_i, \\
  \lsqrt{J_{ij}} &= \tria \left( \begin{bmatrix} A_i^\T \Xi_{22} & \lsqrt{J_i} \end{bmatrix} \right),
\end{align}
\end{subequations}
where \(\Xi_{11}\), \(\Xi_{21}\), \(\Xi_{22}\) are defined via
\begin{equation}
  \begin{bmatrix}
    \Xi_{11} & 0 \\
    \Xi_{21} & \Xi_{22}
  \end{bmatrix}
  = \operatorname{tria} \left(
    \begin{bmatrix}
      \lsqrt{C_i}^\T \lsqrt{J_j} & I \\
      \lsqrt{J_j} & 0
    \end{bmatrix} \right).
\end{equation}
See \citet{Yaghoobi2022} for the detailed derivation.

\paragraph{The filtering marginals.}
By computing a cumulative sum of the elements \(a_{1:N}\) with the binary operation \(\otimes_f\), we obtain the filtering marginals at time \(n=1, \dots, N\) as
\begin{equation}
  p(Y_n | Z_{1:n}) = \mathcal{N} \left(Y_n; m_n^f, P_n^f \right), \qquad \text{with} \qquad
  m_n^f \coloneqq b_{1:n}, \quad
  \lsqrt{P_n^f} \coloneqq  \lsqrt{C_{1:n}}.
\end{equation}
This concludes the \rone{time-parallel} linear Gaussian square-root filter.

\subsubsection{\rone{Time-Parallel} Linear Gaussian Smoothing in Square-Root Form}
Similarly to the filtering equations, the linear Gaussian smoothing can also be formulated in terms of smoothing elements \(b_n\) and an associative operator \(\otimes_s\), and the smoothing marginals can also be computed with a parallel prefix-sum algorithm.

\paragraph{Parameterization of the smoothing elements.}
The smoothing elements \(b_n\) can be described by a set of parameters \(\{E_n, g_n, \lsqrt{L_n}\}\), as
\begin{equation}
  \label{eq:sqrt-smoothing-elements}
  b_n = p(Y_n \mid Z_{1:n}, Y_{n+1}) = \mathcal{N} \left( Y_n; E_n Y_{n+1} + g_n, \lsqrt{L_n} \lsqrt{L_n}^\T \right).
\end{equation}
The smoothing element parameters can be computed as
\begin{subequations}
\begin{align}
  E_n &= \Pi_{21} \Pi_{11}^{-1}, \\
  g_n &= m_n^f - E_n \Phi_n m_n^f, \\
  \lsqrt{L_n} &= \Pi_{22},
\end{align}
\end{subequations}
where \(I\) is the identity matrix and the matrices \(\Pi_{11}\), \(\Pi_{21}\), \(\Pi_{22}\) are defined via
\begin{equation}
  \begin{bmatrix}
    \Pi_{11} & 0 \\
    \Pi_{21} & \Pi_{22}
  \end{bmatrix}
  = \operatorname{tria} \left(
    \begin{bmatrix}
      \Phi_n \lsqrt{P_n^f} & \lsqrt{Q_n} \\
      \lsqrt{P_n^f} & 0
    \end{bmatrix} \right).
\end{equation}

\paragraph{Associative smoothing operator.}
Given two smoothing elements \(b_i, b_j\), parameterized in square-root form by \(b_i = \{E_i, g_i, \lsqrt{L_i}\}\) and \(b_j = \{E_j, g_j, \lsqrt{L_j}\}\), the associative smoothing operator \(\otimes_s\) computes the smoothing element \(b_{ij} = b_i \otimes_s b_j\) as
\begin{subequations}
  \label{eq:sqrt-associative-smoothing}
\begin{align}
  E_{ij} &= E_i E_j, \\
  g_{ij} &= E_i g_j + g_i, \\
  \lsqrt{L_{ij}} &= \tria \left( \begin{bmatrix} E_i \lsqrt{L_j} & \lsqrt{L_i} \end{bmatrix} \right).
\end{align}
\end{subequations}

\paragraph{The smoothing marginals.}
The smoothing marginals can then be retrieved from the reverse cumulative sum of the smoothing elements as
\begin{subequations}
\begin{align}
  p(Y_n \mid Z_{1:N}) &= \mathcal{N} \left( Y_n; m_n^s, P_n^s \right), \\
  m_n^s &= g_{n:N}, \\
  \lsqrt{P_n^s} &= \lsqrt{L_{n:N}}.
\end{align}
\end{subequations}
Refer to \citet{Yaghoobi2022} for a thorough derivation.
This concludes the \rone{time-parallel} Rauch--Tung--Striebel smoother, \rone{which can can be used as a parallel-in-time probabilistic numerical ODE solver for affine ODEs.}
The full algorithm
is summarized in \cref{alg:parallel-rts}.

\begin{algorithm}[H]
  \caption{\rone{Time-parallel} Rauch--Tung--Striebel Smoother (\texttt{\rone{ParaRTS}})}
  \label{alg:parallel-rts}
  \begin{algorithmic}[1]
    \Require
    Initial distribution \((\mu_0, \Sigma_0)\),
    linear transition models \(\{(\Phi_n, Q_n)\}_{n=0}^{N-1}\),
    affine observation models \(\{(H_n, d_n)\}_{n=1}^N\),
    data \(\{Z_n\}_{n=1}^N\).
    \Ensure Smoothing marginals \(p(Y_n \mid Z_{1:N}) = \mathcal{N}\!\left( Y_n; m_n^s, P_n^s \right)\)

    \Function{\texttt{ParaRTS}}{\((\mu_0, \Sigma_0)\), \(\{(\Phi_n, Q_n)\}_{n=1}^N\), \(\{(H_n, d_n)\}_{n=1}^N\), \(\{Z_n\}_{n=1}^N\)}

    \State \emph{Compute the filtering elements:}
    \Statex \(a_n = ( A_n, b_n, \lsqrt{C_n}, \eta_n, \lsqrt{J_n} )\) for all \(n = 1, \dots, N\)
    \Comment{Eq.~(\ref{eq:sqrt-filtering-elements})}

    \State \emph{Run the time-parallel Kalman filter:}
    \Statex \(\left\{ \left( A_n, b_n, \lsqrt{C_n}, \eta_n, \lsqrt{J_n} \right) \right\}_{n=1}^N \leftarrow \texttt{AssociativeScan} \left( \otimes_f, (a_n)_{n=1}^N \right)\)
    \Comment{Eq.~(\ref{eq:sqrt-associative-filtering})}
    \Statex \(p \left( Y_n \mid Z_{1:N} \right) = \mathcal{N}\!\left( Y_n; m_n^f, P_n^f \right) \leftarrow \mathcal{N}\!\left( Y_n; b_n, C_n \right) \)
    \Comment{Filtering marginals}

    \State \emph{Compute the smoothing elements:}
    \Statex \(b_n = ( E_n, g_n, \lsqrt{L_n} )\) for all \(n = 0, \dots, N\)
    \Comment{Eq.~(\ref{eq:sqrt-smoothing-elements})}

    \State \emph{Run the time-parallel Rauch--Tung--Striebel smoother:}
    \Statex \(\left\{ \left (E_n, g_n, \lsqrt{L_n} \right) \right\}_{n=1}^N \leftarrow \texttt{ReverseAssociativeScan} \left( \otimes_s, (b_n)_{n=1}^N \right)\)
    \Comment{Eq.~(\ref{eq:sqrt-associative-smoothing})}
    \Statex \Return \(\mathcal{N}\!\left( Y_n; g_n, L_n \right) \) for all \(n = 1, \dots, N\)
    \Comment{Smoothing marginals}

    \EndFunction

  \end{algorithmic}
\end{algorithm}

\subsection{\rone{Time-Parallel} Approximate Inference in Nonlinear Vector Fields}

Let us now consider the general case: An IVP with nonlinear vector field
\begin{equation}
  \dot{y}(t) = f(y(t), t), \quad t \in [0, T], \qquad y(0) = y_0.
\end{equation}
As established in \cref{sec:background}, the corresponding state estimation problem is
\begin{subequations}
  \label{eq:nonlinear-estimation}
  \begin{align}
    Y_0 &\sim \mathcal{N} \left( \mu_0, \Sigma_0 \right), \\
    Y_{n+1} \mid Y_n &\sim \mathcal{N} \left( \Phi_n Y_n, Q_n \right), \\
    Z_n \mid Y_n &\sim \delta \left( E_1 Y_n - f \left( E_0 Y_n, t_n \right) \right), \label{eq:nonlinear-estimation-obs}
  \end{align}
\end{subequations}
with temporal discretization \(\mathbb{T} = \{t_n\}_{n=1}^N \subset [0, T]\) and
zero data \(Z_n = 0\) for all \(n = 1, \dots, N\).
In this section, we describe a parallel-in-time algorithm for solving this state estimation problem:
the \emph{\rone{time-parallel} iterated extended Kalman smoother}.

\subsubsection{Globally Linearizing the State-Space Model}
\label{sec:linearizing}
To make inference tractable, we will linearize the whole state-space model along a reference trajectory.
And since the observation model (specified in \cref{eq:nonlinear-estimation-obs}) is the only nonlinear part of the state-space model, it is the only part that requires linearization.
In this paper, we only consider linearization with a first-order Taylor expansion, but other methods are possible; see
\cref{remark:approx-linearization,remark:statistical-linearization}.

For any time-point \(t_n \in \mathbb{T}\), we approximate the nonlinear observation model
\rone{\cref{eq:nonlinear-estimation-obs}}
with an affine observation model by performing a first-order Taylor series expansion around a linearization point \(\eta_n \in \mathbb{R}^{d (\nu+1)}\). We obtain the affine model%
\begin{equation}
  \label{eq:linearized-observation-model}
    Z_n \mid Y_n \sim \delta \left( H_n Y_n - d_n \right),
\end{equation}
with \(H_n\) and \(d_n\) defined as
\begin{subequations}
  \begin{align}
    H_n &\coloneqq E_1 - F_y(E_0 \eta_n, t_n) E_0, \\
    d_n &\coloneqq f(E_0 \eta_n, t_n) - F_y(E_0 \eta_n, t_n) E_0 \eta_n,
  \end{align}
\end{subequations}
where \(F_y\) denotes the Jacobian of \(f\) with respect to \(y\).

In iterated extended Kalman smoothing, this linearization is performed \emph{globally} on all time steps simultaneously along a trajectory of linearization points \( \left\{ \eta_n \right\}_{n=1}^N \subset \mathbb{R}^{d(\nu+1)} \).
We obtain the following linearized inference problem:%
\begin{subequations}
  \label{eq:linearized-inference-problem}
  \begin{align}
    Y_0 &\sim \mathcal{N} \left( \mu_0, \Sigma_0 \right), \\
    Y_{n+1} \mid Y_n &\sim \mathcal{N} \left( \Phi_n Y_n, Q_n \right), \\
    Z_n \mid Y_n &\sim \delta \left( H_n Y_n - d_n \right),
  \end{align}
\end{subequations}
with zero data \(Z_n = 0\) for all \(n = 1, \dots, N\).
This is now a linear state-space model with linear Gaussian observations.
It can therefore be solved exactly with the numerically stable, time-parallel Kalman filter and smoother presented in \cref{sec:affine}.

\begin{remark}[Linearizing with approximate Jacobians (\texttt{EK0} \& \texttt{DiagonalEK1})]
  \label{remark:approx-linearization}
  To reduce the computational complexity with respect to the state dimension of the ODE,
  the vector field can also be linearized with an approximate Jacobian.
  Established choices include \(F_y \approx 0\)
  and \(F_y \approx \operatorname{diag}(\nabla_y f)\),
  which result in probabilistic ODE solvers known as the \texttt{EK0} and \texttt{DiagonalEK1}, respectively.
  See \citet{Kramer2022highdim} for more details.
\end{remark}

\begin{remark}[Statistical linear regression]
  \label{remark:statistical-linearization}
  Statistical linear regression (SLR) is a more general framework for approximating conditional distributions with affine Gaussian distributions, and many well-established filters can be understood as special cases of SLR.
  This includes notably the Taylor series expansion used in the EKF/EKS, but also sigma-point methods such as the unscented Kalman filter and smoother
  \citep{julier2000ufk1,julier2004ufk2,sarkka2008ufs}, and more.
  For more information on SLR-based filters and smoothers refer to
  \citet[Chapter 9]{sarkka2023bfs2}.
\end{remark}

\subsubsection{Iterated Extended Kalman Smoothing}
\label{sec:ieks}

The IEKS
\citep{Bell1994,sarkka2023bfs2}
is an approximate Gaussian inference method for nonlinear state-space models, which iterates between linearizing the state-space model along the current best-guess trajectory and computing a new state trajectory estimate by solving the linearized model exactly.
It can equivalently also be seen as an efficient implementation of the Gauss--Newton method, applied to maximizing the posterior density of the state trajectory
\citep{Bell1994}.
This also implies that the IEKS computes not just some Gaussian estimate, but the \emph{maximum a posteriori} (MAP) estimate of the state trajectory.
In the context of probabilistic numerical ODE solvers, the IEKS has been previously explored by \citet{Tronarp2021a}, and the resulting MAP estimate has been shown to satisfy polynomial convergence rates to the true ODE solution.
Here, we formulate an IEKS-based probabilistic ODE solver in a parallel-in-time manner, by exploiting the time-parallel formulation of the Kalman filter and smoother from \cref{sec:affine}.

The IEKS is an iterative algorithm, which starts with an initial guess of the state trajectory and then iterates between the following two steps:
\begin{enumerate}
\item \emph{Linearization step:} Linearize the state-space model along the current best-guess trajectory. This can be done independently for each time step and is therefore fully parallelizable.
\item \emph{Linear smoothing step:} Solve the resulting linear state-space model exactly with the time-parallel Kalman filter and smoother from \cref{sec:affine}.
\end{enumerate}
The algorithm terminates when a stopping criterion is met, for example when the change in the MAP estimate between two iterations is sufficiently small.
A pseudo-code summary of the method is provided in \cref{alg:solver}.

\begin{algorithm}
  \caption{Parallel-in-Time Probabilistic Numerical ODE Solver \rone{(\texttt{ParaIEKS})}}
  \label{alg:solver}
  \begin{algorithmic}[1]
    \Require
    ODE-IVP \((f, y_0)\),
    prior transition model \(\left\{(\Phi_n, Q_n)\right\}_{n=1}^N\),
    time grid \(\left\{ t_n \right\}_{n=0}^N \subset [0, T]\).

    \Ensure Probabilistic numerical ODE solution \\
    \(p \left( y(t_n) ~\Big|~ y(0) = y_0, \left\{ \dot{y}(t_m) = f(y(t_m), t_m) \right\}_{m=1}^N \right)
    \approx \mathcal{N} \left( \mu_n, \Sigma_n \right)\) for \(n = 1, \dots, N\).

    \Function{\texttt{ParaIEKS}}{\((f, y_0)\), \(\left\{(\Phi_n, Q_n)\right\}_{n=1}^N\), \(\left\{ t_n \right\}_{n=0}^N\)}

    \State \(\mu_0,\ \Sigma_0 \leftarrow \texttt{ComputeExactInitialState}(f, y_0),\ 0\)
    \Comment{With automatic differentiation}
    \State \(\eta_n \leftarrow \mu_0\) for all \(n = 0, \dots, N\)
    \Comment{Constant initial guess of state trajectory}
    \While {stopping criterion not met}
      \LComment{Linearize the state-space model as in Sec.~\ref{sec:linearizing} (fully parallelizable):}
      \For{\(n = 1, \dots, N\)}
        \State
        \( H_n, d_n \leftarrow \texttt{LinearizeObservationModel}(f, \eta_n, t_n) \)
      \EndFor
      \LComment{Run the time-parallel RTS in the linearized model as in Sec.~\ref{sec:affine}:}
      \State
      \( \left\{\mu_n, \Sigma_n \right\}_{n=1}^N \leftarrow \texttt{\rone{ParaRTS}} \left( (\mu_0, \Sigma_0), \{(\Phi_n, Q_n)\}_{n=1}^N, \{(H_n, d_n)\}_{n=1}^{N}, \{0\}_{n=1}^N \right) \)
      \LComment{Choose the mean trajectory as the next linearization trajectory:}
      \State \(\eta_n \leftarrow \mu_n\) for all \(n = 1, \dots, N\)
    \EndWhile

    \State \Return \( y(t_n) \sim \mathcal{N} \left( E_0 \mu_n, E_0 \Sigma_n E_0^\T \right)\) for \(n = 1, \dots, N\).

    \EndFunction
  \end{algorithmic}
\end{algorithm}

As with the sequential filtering-based probabilistic ODE solvers as presented in \cref{sec:background}, the mean and covariance of the initial distribution \(Y(0) \sim \mathcal{N}\left( \mu_0, \Sigma_0 \right)\) are chosen such that \(\mu_0\) corresponds to the exact solution of the ODE and its derivatives and \(\Sigma_0\) is set to zero; see also \citet{Kramer2020stable}.
The initial state trajectory estimate \(\{\eta_n\}_{n=0}^N\) is chosen to be constant, that is, \(\eta_n = \mu_0\) for all \(n = 0, \dots, N\).
Note that since only \(E_0 \eta_n\) is required to perform the linearization, it could equivalently be set to \(\eta_n = [y_0, 0, \dots, 0]\) for all \(n\).

\begin{remark}[Other choices for the initial state trajectory]
  \label{remark:initialization}
  In the sequential IEKS the initial state trajectory is typically computed with a standard sequential extended Rauch--Tung--Striebel smoother \citep{Tronarp2021a}.
  But in the context of the time-parallel IEKS, we can not afford to apply a sequential method first as this would break the parallel-in-time nature of the algorithm and increase the computational complexity from logistic back to linear in the number of time points.
  \rone{
  One possible alternative choice is compute an initial state trajectory with a sequential method on a coarse time grid with \(\order{\log(N)}\) time points, similar to classic parallel-in-time methods such as Parareal \citep{LIONS2001}.
  But, this introduces additional hyperparameters and it is only beneficial if this trajectory is indeed better than a constant guess, which often requires a certain minimum number of time points.}
  We therefore choose the constant initial guess for the state trajectory in this paper, and leave the exploration of more sophisticated initialization strategies for future work.
\end{remark}

Finally, the stopping criterion should be chosen such that the algorithm terminates when the MAP estimate of the state trajectory has converged.
In our experiments, we chose a combination of two criteria:
(i) the change in the state trajectory estimate between two iterations is sufficiently small, or
(ii) the change in the \emph{objective value} between two iterations is sufficiently small, where the objective value is defined as the negative log-density of the state trajectory under the prior:
\begin{equation}
  \label{eq:objective-value}
  \mathcal{V}(\eta_{0:N}) = \frac{1}{2} \sum_{n=1}^N \left\| \eta_n - \Phi_n \eta_{n-1} \right\|_{Q^{-1}_n}^2.
\end{equation}
In our experiments, we use a relative tolerance of \(10^{-13}\) for the first criterion and absolute and relative tolerances of \(10^{-9}\) and \(10^{-6}\) for the second criterion, respectively.
Note that this procedure is compatible with the uncertainty calibration mentioned in \cref{sec:practical-considerations} as the mean estimator is agnostic to the value of \(\sigma\) (and therefore \(Q\)).

\rtwo{
\begin{remark}[Convergence of the time-parallel IEKS]
The IEKS is equivalent to the Gauss--Newton method for computing the MAP estimate of the state trajectory
\citep{Bell1994},
and the MAP estimate has been shown to satisfy polynomial convergence rates to the true ODE solution,
\citep{Tronarp2021a}.
And under mild conditions on the Jacobian of the vector field, the MAP is a local optimum and the Gauss--Newton method is locally convergent
\citep{Tronarp2021a,Knoth1989}.
While therefore the initial guess of the state trajectory must be sufficiently close to the true solution for the method to be guaranteed to converge, in our experiments we found that the constant initialization works well in practice; see \cref{sec:experiments} below.
In addition, the stability of the method could be further improved by using a more sophisticated optimizer, such as the Levenberg--Marquardt method or Gauss--Newton with line search \citep{sarkka_lm_linesearch}, or the alternating method of multipliers
\citep[ADMM;][]{Gao2019}.
We leave the exploration of these methods for future work.
\end{remark}
}

\subsection{\rone{Computational Complexity}}
\label{sec:complexity}

Let
\(C_\text{KS-step}^s\)
be the summed costs of a predict, update, and smoothing step in the sequential Kalman smoother, and
let
\(C_\text{KS-step}^p\)
be the corresponding cost in the time-parallel formulation
(which differs by a constant factor, i.e.\ \(C_\text{KS-step}^p \propto C_\text{KS-step}^s\)).
Let \(C_\text{linearize}\) be the cost of linearizing the vector field at a single time point, which requires both evaluating the vector field and computing its Jacobian.
Then, on a grid with \(N\) time points, the cost of a standard sequential extended Kalman smoother is linear in \(N\), with
\begin{equation}
  C_\text{EKS}^s = N \cdot C_\text{linearize} + N \cdot C_\text{KS-step}^s.
\end{equation}
The computational cost of the time-parallel IEKS differs in two ways:
(i) the prefix-sum formulation of the Kalman smoother enables a time-parallel inference with logarithmic complexity, and
(ii) the linearization is not done locally in a sequential manner but can be performed globally, fully in parallel.
In the following we assume that we have at least \(N\) processors / threads / cores available (simply referred to as ``cores'' hereafter), which allow us to take perfect advantage of the temporal parallelization,
and we disregard the effect of the number of cores on all other parts of the algorithm (such as the influence of the number of cores on the computation time of matrix multiplications).
The computational cost of the time-parallel IEKS run for \(k\) iterations is then
\begin{equation}
  C_\text{IEKS}^p
  = k \cdot \left( C_\text{linearize} + 2 \log_2(N) \cdot \left( C_\text{KS-step}^p \right) \right).
\end{equation}
Thus, the runtime of the time-parallel IEKS is logarithmic in the number of time points, provided a sufficient number of cores are available.
To make the comparison to the sequential EKS more explicit, we consider the \emph{speedup} of the time-parallel IEKS over the sequential EKS, defined as
\(S = C_\text{EKS}^s / C_\text{IEKS}^p\).
By simply inserting the computational costs from above and re-arranging the terms, it can be shown that the speedup is bounded by
\begin{align}
  S &\leq \min \left\{
      \frac{N}{k} \cdot
      \frac{C_\text{linearize} + C_\text{KS-step}^s}{C_\text{linearize}},
      \frac{N}{2 k \log_2(N)} \cdot
      \frac{C_\text{linearize} + C_\text{KS-step}^s}{C_\text{KS-step}^p}
      \right\}.
\end{align}
This term highlights two different regimes:
If the cost of linearizing the vector field is small compared to the cost of the Kalman smoother,
then the speedup is approximately bounded by
\(S \lesssim \frac{N}{2 k \log_2\!N}\)
(since the sequential and parallel Kalman smoother steps have similar cost).
On the other hand, if the cost of linearizing the vector field is large compared to the cost of the Kalman smoother, then the benefit of parallelizing the linearization step dominates and the speedup is bounded by
\(S \lesssim \frac{N}{k} \).

\section{Experiments}
\label{sec:experiments}

This section investigates the utility and performance of the proposed parallel-in-time ODE filter
on a range of experiments.
It is structured as follows:
First, \cref{sec:parallel-vs-sequential} compares the performance of the parallel-in-time probabilistic numerical ODE solver to its sequential counterpart on multiple test problems.
\rone{
\Cref{sec:parallel-scaling} then investigates the parallel scaling of the proposed method and evaluates it on a range of problem sizes and for different GPUs.
}
Finally, \cref{sec:benchmark} benchmarks the proposed method against other well-established ODE solvers, including both classic and probabilistic numerical methods.

\rone{
In the following, we refer to all filtering-based probabilistic numerical ODE solvers by their underlying filtering and smoothing algorithms:
the proposed parallel-in-time method is referred to as \texttt{ParaIEKS},
the corresponding sequential method is simply \texttt{IEKS},
and the established sequential EKS-based method is referred to as \texttt{EKS}.
And, whenever relevant, we write the smoothness of the chosen IWP prior into the method name, for example as \texttt{ParaIEKS($\nu\!=\!1$)} or \texttt{EKS($\nu\!=\!2$)}.
}

\paragraph{Implementation}
\label{sec:implementation}
All experiments are implemented in the Python programming language
with the JAX software framework
\citep{jax2018github}.
Reference solutions are computed with
SciPy \citep{scipy} and
Diffrax
\citep{kidger2021}.
Unless specified otherwise, experiments are run on an NVIDIA V100 GPU\@.
Code for the implementation and experiments is publicly available on GitHub.%
\footnote{\url{https://github.com/nathanaelbosch/parallel-in-time-ode-filters}}

\subsection{The \rone{Parallel-in-time Solver} Compared to its Sequential Version}
\label{sec:parallel-vs-sequential}

\begin{figure}[t]
  \includegraphics{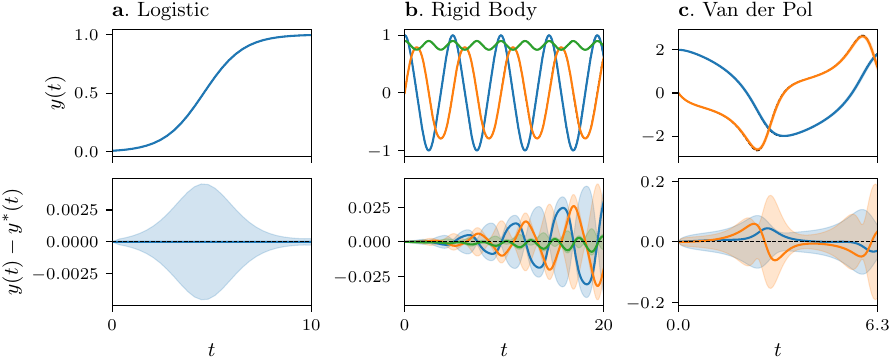}
  \caption{
    \emph{Trajectories, errors, and error estimates computed by the parallel-in-time solver.}
    Top row: ODE solution trajectories.
    Visually, all three test problems seem to be solved accurately.
    Bottom row: Numerical errors (lines) and error estimates (shaded area).
    Ideally, for good calibration, the error should be of similar magnitude than the error estimate.
    The posterior appears underconfident on the logistic equation, and well-calibrated on the rigid body and Van der Pol problems.
  }
  \label{fig:ode-solutions}
\end{figure}

We first compare the proposed parallel-in-time ODE solver \rone{\texttt{ParaIEKS}} to its sequential counterpart: a probabilistic solver based on the sequential implementation \rone{\texttt{IEKS}}.
We evaluate the solvers on three test problems:
The logistic ordinary differential equation
\begin{equation}
  \label{eq:logistic}
  \dot{y}(t) = y(t) \left( 1 - y(t) \right), \qquad t \in [0, 10], \qquad
  y(0) = 0.01,
\end{equation}
an initial value problem based on the rigid body dynamics
\citep{hairer2008solving}
\begin{equation}
  \dot{y}(t) = \begin{bmatrix} -2 y_2(t) y_3(t) \\ 1.25 y_1(t) y_3(t) \\ -0.5 y_1(t) y_2(t) \end{bmatrix}, \qquad t \in [0, 20], \qquad
  y(0) = \begin{bmatrix} 1 \\ 0 \\ 0.9 \end{bmatrix},
\end{equation}
and the Van der Pol oscillator \citep{van1920}
\begin{equation}
  \dot{y}(t) = \begin{bmatrix} y_2(t) \\ \mu \left( \left( 1 - y_1(t)^2 \right) y_2(t) - y_1(t) \right) \end{bmatrix}, \qquad t \in [0, 6.3], \qquad
  y(0) = \begin{bmatrix} 2 \\ 0 \end{bmatrix},
\end{equation}
here in a non-stiff version with parameter \(\mu = 1\).

We first solve the three problems with \rone{\texttt{ParaIEKS}} on grids of sizes \(30\), \(150\), and \(100\), respectively for the logistic, rigid body, and Van der Pol problem, with a two-times integrated Wiener process prior.
Reference solutions are computed with diffrax's \texttt{Kvaerno5} solver using adaptive steps and very low tolerances \(\tau_{\{\text{abs},\text{rel}\}}=10^{-12}\)
\citep{kidger2021, kvaerno2004}.
\Cref{fig:ode-solutions} shows the resulting solution trajectories, together with numerical errors and error estimates.
For these grid sizes, \rone{\texttt{ParaIEKS}} computes accurate solutions on all three problems.
Regarding calibration, the posterior appears underconfident for the logistic equation, as it overestimates the numerical error by more than one order of magnitude, but is reasonably confident on the on the rigid body and Van der Pol problems where the error estimate is of similar magnitude as the numerical error.

\begin{figure}[p]
  \includegraphics[width=\linewidth]{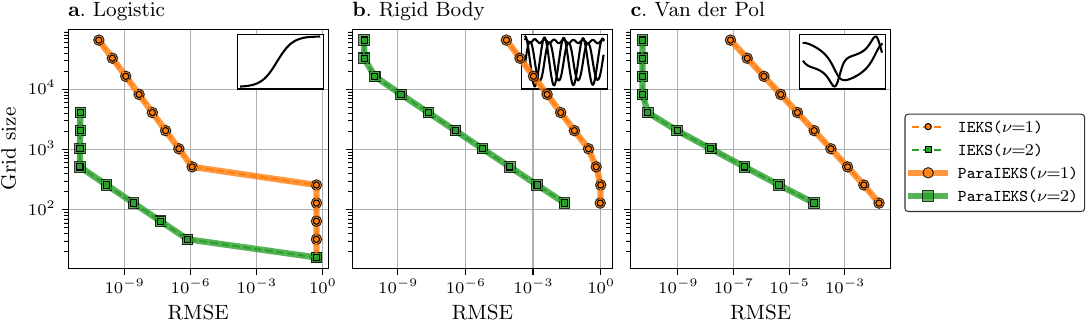}
  \caption{
    \emph{The sequential and parallel IEKS compute numerically identical solutions.}
    For all three problems and all considered grid sizes, the sequential \rone{\texttt{IEKS}} and the parallel \rone{\texttt{ParaIEKS}} achieve (numerically) identical errors.
    This is expected, as both versions compute the same quantities and only differ in their implementation.
  }
  \label{fig:seq-vs-par-ns}
\end{figure}

\begin{figure}[p]
  \includegraphics[width=\linewidth]{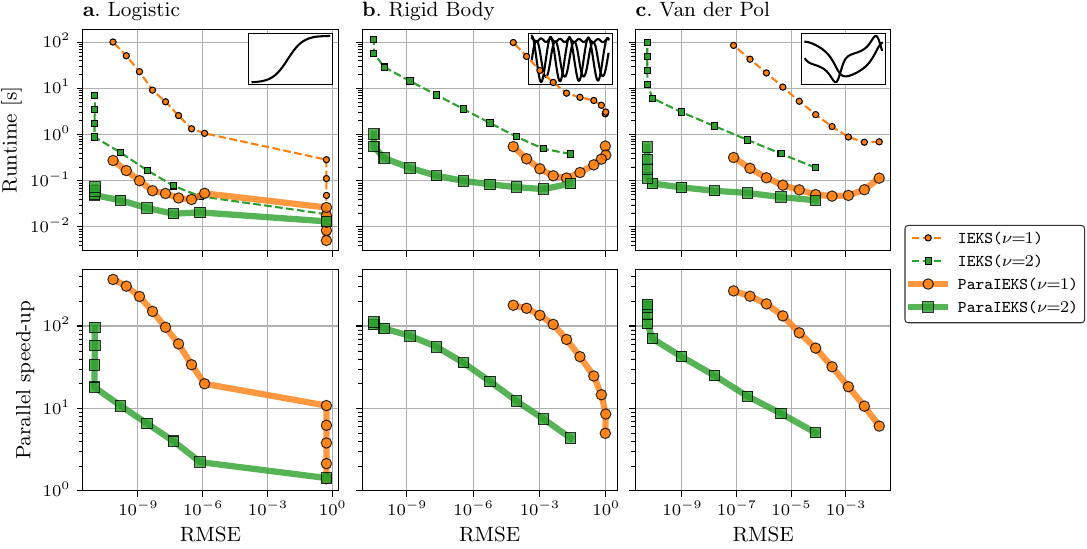}
  \caption{
    \emph{Work-precision diagrams for the sequential and parallel IEKS-based ODE solver.}
    Top row:
    Runtime in seconds per error (lower-left is better).
    Bottom row: Speed-up of the parallel \rone{\texttt{ParaIEKS}} over the sequential \rone{\texttt{IEKS}} (higher is better).
    Across all problems, grid sizes, and priors, \rone{\texttt{ParaIEKS}} outperforms the sequential \rone{\texttt{IEKS}}.
  }
  \label{fig:seq-vs-par}
\end{figure}

Next, we investigate the performance of \rone{\texttt{ParaIEKS}} and compare it to its sequential implementation \rone{\texttt{IEKS}}.
We solve the three problems with both methods on a range of grid sizes, with both a one- and two-times integrated Wiener process prior.
Reference solutions are computed with diffrax's \texttt{Kvaerno5} solver using adaptive steps and very low tolerances
(\(\tau_{\text{abs}}=10^{-16}, \tau_{\text{rel}}=10^{-13}\)).
\Cref{fig:seq-vs-par-ns} shows the achieved root-mean-square errors (RMSE) for different grid sizes in a work-precision diagram.
As expected, \rone{\texttt{ParaIEKS} and \texttt{IEKS}} always achieve the same error for each problem and grid size, as both versions use the same initialization, compute the same quantities, and only differ in their filter and smoother implementation.
However, the methods differ significantly in actual runtime, as shown in \cref{fig:seq-vs-par}.
In our experiments on an NVIDIA V100 GPU, \rone{\texttt{ParaIEKS}} is always strictly faster than \rone{\texttt{IEKS}} across all problems, grid sizes, and priors,
and we observe speed-ups of multiple orders of magnitude.
Thus, when working with a GPU, \rone{\texttt{ParaIEKS}} appears to be strictly superior to the sequential \rone{\texttt{IEKS}}.

\subsection{\rone{Runtimes and Parallel Scaling of the Parallel-in-Time Solver}}
\label{sec:parallel-scaling}
\begin{figure}[t]
  \includegraphics{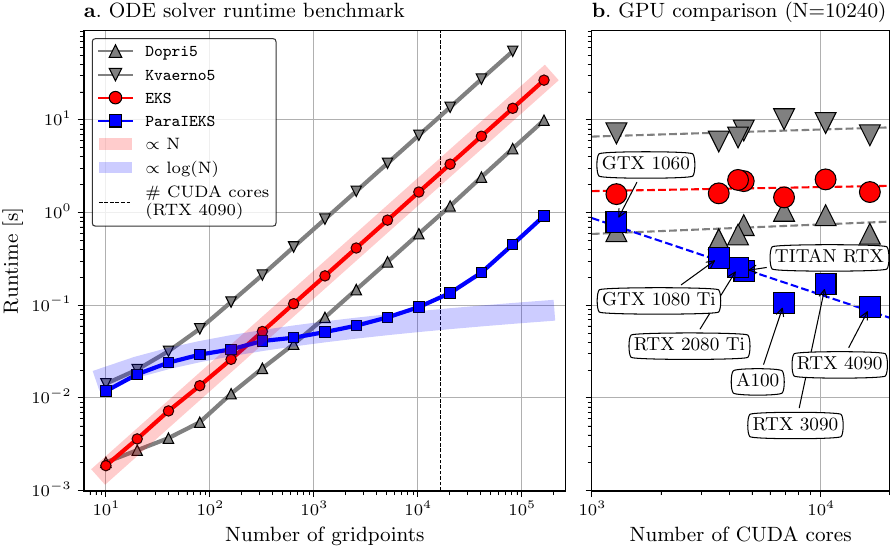}
  \caption{
    \emph{The time-parallel \texttt{ParaIEKS} shows logarithmic scaling in the grid size and benefits from GPU improvements.}
    In comparison, the sequential \texttt{IEKS} and the classic \texttt{Dopri5} and \texttt{Kvaerno5} solvers show the expected linear runtime complexity (left).
    The sequential methods also do not show relevant changes in runtime for GPUs with more CUDA cores, whereas \texttt{ParaIEKS}'s runtime improves (right).
  }
  \label{fig:parallel-scaling}
\end{figure}

Next we evaluate how the runtime of the proposed method scales with respect to increasing problem sizes, as well as increased numbers of available CUDA cores.
We compare \rone{\texttt{ParaIEKS}} to the established sequential \rone{\texttt{EKS}}, as well as to two classic ODE solvers:
the explicit Runge--Kutta method \texttt{Dopri5} \citep{dormand1980,shampine1986},
and the implicit Runge--Kutta method \texttt{Kvaerno5} \citep{kvaerno2004};
both provided by the Diffrax python package \citep{kidger2021}.
We apply all methods to the logistic ODE on a range of grid sizes,
resulting from time discretizations with step sizes
\rone{\(h = 2^0, 2^{-1}, \dots, 2^{-14}\)},
as well as for multiple GPUs with varying numbers of CUDA cores.
\Cref{fig:parallel-scaling} shows the results.

First, we observe the expected logarithmic scaling of \rone{ParaIEKS} with respect to the grid size, for grids of size up to around \(5 \cdot 10^3\)
(\cref{fig:parallel-scaling}a).
For larger grid sizes the runtime of \rone{\texttt{ParaIEKS}} starts to grow linearly.
This behavior is expected: The NVIDIA RTX 4090 used in this experiment has 16384 CUDA cores, so for larger grids the filter and smoother pass can not be fully parallelized anymore and additional grid points need to be processed sequentially.
Note also that for very large grids, the \(24\text{GB}\) memory of the GPU can become a limiting factor.
Nevertheless, for large grid sizes \(N > 10^3\), \rone{\texttt{ParaIEKS}} achieves the lowest runtimes out of all the tested solvers, with up to an order of magnitude faster ODE solutions.

\Cref{fig:parallel-scaling}b.\ shows runtimes for different GPUs with varying numbers of CUDA cores for a grid of size \(N = 10240\).
We observe that the sequential \rone{\texttt{IEKS}}, \texttt{Dopri5} and \texttt{Kvaerno5} solvers do not show a benefit from the improved GPU hardware, which is expected as the method does not leverage parallelization.
On the other hand, the runtime of \rone{\texttt{ParaIEKS}} decreases as the number of CUDA cores increases, and we observe speed-ups of up to an order of magnitude by using a newer GPU with a larger number of CUDA cores.

\subsection{Benchmarking the \rone{Parallel-in-time Probabilistic Numerical ODE Solver}}
\label{sec:benchmark}
Finally, we compare \rone{\texttt{ParaIEKS}} to a range of well-established \rtwo{implicit} ODE solvers, including both classic and probabilistic numerical methods:
we compare against the implicit Euler method (\texttt{ImplicitEuler}) and the \texttt{Kvaerno3} and \texttt {Kvaerno5} solvers \citep{kvaerno2004} provided by Diffrax \citep{kidger2021},
as well as the sequential \rone{\texttt{EKS}} with local linearization, which is one of the currently most popular probabilistic numerical ODE solvers.
We evaluate all solvers on all three test problems, namely the logistic ODE, the rigid body problem, and the Van der Pol oscillator, as introduced in \cref{sec:parallel-vs-sequential}.
Reference solutions are computed with diffrax's \texttt{Kvaerno5} solver with adaptive steps and a very low error tolerance setting (\(\tau_{\text{abs}}=10^{-16}, \tau_{\text{rel}}=10^{-13}\)).

\begin{figure}[t]
  \includegraphics[width=1.0\linewidth]{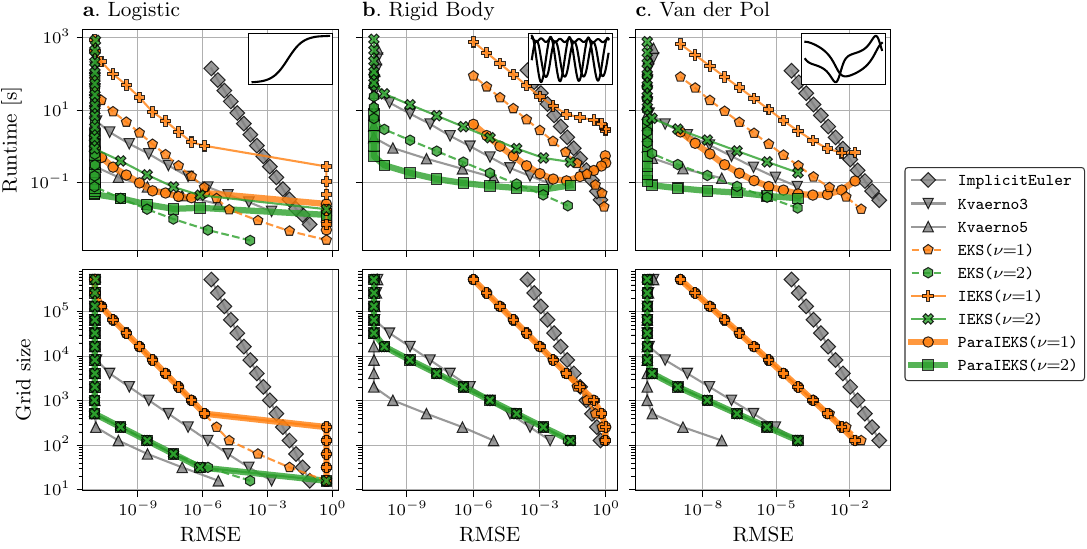}
  \caption{
    \emph{Benchmarking \rone{\texttt{ParaIEKS}} against other commonly used numerical ODE solvers.}
    Top row: Work-precision diagrams showing runtimes per error for a range of different ODE solvers (lower-left is better).
    Bottom row: Errors per specified grid size (lower-left is better).
    Per grid size, the closely related \rone{\texttt{EKS}, \texttt{IEKS}, and \texttt{ParaIEKS}} solvers often coincide; \texttt{Kvaerno5} achieves the lowest error per step as it has the highest order.
    In terms of runtime, \rone{\texttt{ParaIEKS}} outperforms all other methods on medium-to-high accuracy settings due to its logarithmic time complexity.
  }
  \label{fig:wpd}
\end{figure}

\begin{figure}[t]
  \includegraphics[width=1.0\linewidth]{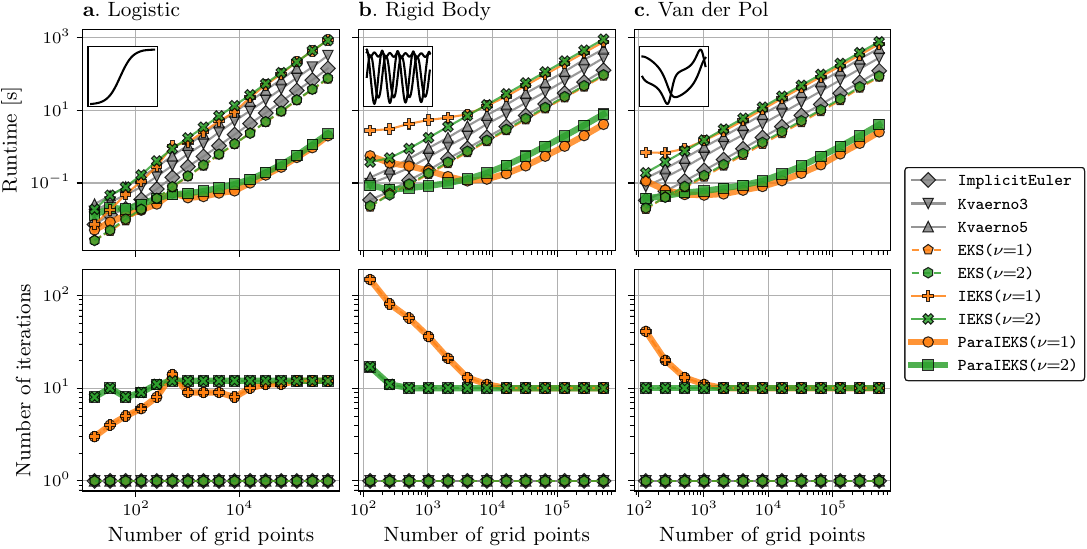}
  \caption{
    \emph{Runtimes of the ODE solvers for each grid size, and number of IEKS iterations.}
    While all sequential solvers demonstrate linear scaling with the number of grid points, \rone{\texttt{ParaIEKS}} shows sub-linear scaling up to a certain grid size (top).
    The number of IEKS iterations until convergence can vary with the grid size and the problem, but it seems that in many cases around \(10\) iterations suffice (bottom).
    The sequential methods solve the ODE in one sweep.
  }
  \label{fig:scaling}
\end{figure}

\Cref{fig:wpd} shows the results as work-precision diagrams.
For small grid sizes (low accuracy), the logarithmic time complexity of \rone{\texttt{ParaIEKS}} seems to not be very relevant and it is outperformed by the non-iterated \rone{\texttt{EKS}}.
In the particular case of the logistic ODE, it further seems that on coarse grids the MAP estimate computed by \rone{\texttt{ParaIEKS}} (and by the sequential \texttt{IEKS}) differs significantly from the actual ODE solution, and thus the error of \rone{\texttt{ParaIEKS}} on coarse grids is high (lower left figure).
However, for larger grid sizes (medium-to-high accuracy), \rone{\texttt{ParaIEKS}} outperforms both its sequential, non-iterated counterpart \rone{\texttt{EKS}}, as well as the classic implicit methods.
In particular, \rone{\texttt{ParaIEKS}} with IWP(2) prior often shows runtimes lower than those of the classic \texttt{Kvaerno5} method for comparable errors, even though it has a lower order of convergence and it requires multiple iterations per solve.
\Cref{fig:scaling} shows the runtimes per grid size of all methods, as well as the number of iterations performed by the IEKS.
We again observe the initial logarithmic scaling of the runtime with increasing problem size, up to a threshold where the runtime starts scaling linearly, as previously shown in \cref{sec:parallel-scaling}.
Overall, the logarithmic time complexity of the proposed \rone{\texttt{ParaIEKS}} appears to be very beneficial for high accuracy settings on GPUs and makes \rone{\texttt{ParaIEKS}} a very competitive ODE solver in this comparison.

\section{Conclusion}
\label{sec:conclusion}

In this work, we have developed a \emph{parallel-in-time} probabilistic numerical ODE solver.
The method builds on iterated extended Kalman smoothing to compute the maximum a posteriori estimate of the probabilistic ODE solution,
and by using the time-parallel formulation of the IEKS it is able to efficiently leverage modern parallel computer hardware such as GPUs to parallelize its computations.
Given enough processors or cores, the proposed algorithm shares the logarithmic cost per time step of the parallel IEKS and the underlying parallel prefix-sum algorithm, as opposed to the linear time complexity of standard, sequentially-operating ODE solvers.
We evaluated the performance of the proposed method in a number of experiments,
and have seen that the proposed parallel-in-time solver can provide speed-ups of multiple orders of magnitude over the sequential IEKS-based solver.
We also compared the proposed method to a range of well-established, both probabilistic and classical ODE solvers, and we have shown that the proposed parallel-in-time method is competitive with respect to the state-of-the-art in both accuracy and runtime.

This work opens up a number of interesting avenues for future research in the intersection of probabilistic numerics and parallel-in-time methods.
Potential opportunities for improvement include the investigation of other optimization algorithms, such as Levenberg--Marquart or ADMM, or the usage of line search, all of which have been previously proposed for the sequential IEKS
\citep{sarkka_lm_linesearch,Gao2019}.
Furthermore, combining the solver with adaptive grid refinement approaches could also significantly improve its performance in practice.
A different avenue would be to extend the proposed method to other related differential equation problems for which sequentially-operating probabilistic numerical methods already exist, such as higher-order ODEs, differential-algebraic equations, or boundary value problems
\citep{Bosch2022,kraemer2021bvp},
or to apply the method to other types of probabilistic inference problems, such as parameter or latent-force inference in ODEs
\citep{Schmidt2021,2022fenrir,beck2024diffusion}.
Finally, the improved utilization of GPUs by our parallel-in-time method could be particularly beneficial to applications in the field of machine learning,
where GPUs are often required to accelerate the computations of deep neural networks.
In summary, the proposed parallel-in-time probabilistic numerical ODE solver not only advances the efficiency of probabilistic numerical ODE solvers, but also paves the way for a range of future research on parallel-in-time probabilistic numerical methods and their application across various scientific domains.

\acks{%
  The authors gratefully acknowledge co-funding by the European Union (ERC, ANUBIS, 101123955. Views and opinions expressed are however those of the author(s) only and do not necessarily reflect those of the European Union or the European Research Council. Neither the European Union nor the granting authority can be held responsible for them). Philipp Hennig is a member of the Machine Learning Cluster of Excellence, funded by the Deutsche Forschungsgemeinschaft (DFG, German Research Foundation) under Germany's Excellence Strategy – EXC number 2064/1 – Project number 390727645; he also gratefully acknowledges the German Federal Ministry of Education and Research (BMBF) through the Tübingen AI Center (FKZ: 01IS18039A); and funds from the Ministry of Science, Research and Arts of the State of Baden-Württemberg.
  The authors would like to thank Research Council of Finland for funding.
  Filip Tronarp was partially supported by the Wallenberg AI, Autonomous Systems and Software Program (WASP) funded by the Knut and Alice Wallenberg Foundation.
  The authors thank the International Max Planck Research School for Intelligent Systems (IMPRS-IS) for supporting Nathanael Bosch.
  The authors are grateful to Nicholas Krämer for many valuable discussion and to Jonathan Schmidt for feedback on the manuscript.
}

{\vskip 0.3in\noindent{\large\bf Individual Contributions}\vskip 0.2in%
  \noindent%
  The original idea for this article came independently from SS and from discussions between FT and NB.
  The joint project was initiated and coordinated by SS and PH.
  The methodology was developed by NB in collaboration with AC, FT, PH, and SS.
  The implementation is primarily due to NB, with help from AC.
  The experimental evaluation was done by NB with support from FT and PH.
  The first version of the article was written by NB, after which all authors reviewed the manuscript.
}

\bibliography{references}

\end{document}